\documentclass[12]{article}
\usepackage{amsfonts}
\usepackage{amsmath}
\usepackage[unicode,colorlinks,plainpages=false]{hyperref}
\usepackage{theorem}

\newtheorem{theorem}{Theorem}

\newcommand{\D}{\mathcal{D}}

\newcommand{\openbox}{$\begin{array}{c}
\hspace*{-0.55em}\sqcap \hspace*{-0.60em}\\[-0.4em] \hline
\multicolumn{1}{c}{\hspace*{-0.60em}}\\[-0.8em]
\end{array}
$}

\usepackage{enumerate}

\begin{document}

\centerline{{\bf On Congruences on Ultraproducts of Algebraic Structures} }

\bigskip
\centerline{Attila Nagy}
\medskip

\begin{abstract}
Let $I$ be a non-empty set and $\mathcal{D}$ an ultrafilter over $I$. For similar algebraic structures $B_i$, $i\in I$ let $\Pi (B_i|i\in I)$ and $\Pi _{\mathcal{D}}(B_i|i\in I)$ denote the direct product and the ultraproduct of $B_i$, respectively. Let $\mathcal{D}^*$ denote the ultraproduct congruence on $\Pi (B_i|i\in I)$. Let the $\wedge$-semilattice of all congruences on an algebraic structure $B$ denoted by ${\bf Con}(B)$.

In this paper we show that, for any similar algebraic structures $A_i$, $i\in I$, there is an embedding $\Phi$ of $\Pi _{\mathcal{D}}({\bf Con}(A_i)|i\in I)$ into ${\bf Con}(\Pi _{\mathcal{D}}(A_i|i\in I)$. We also show that, for every $\sigma \in \Pi ({\bf Con}(A_i)|i\in I)$, the factor algebra $\Pi _{\mathcal{D}}(A_i|i\in I)/\Phi (\sigma /\mathcal{D}^*)$ is isomorphic to $\Pi _{\mathcal{D}}(A_i/\sigma (i)|i\in I)$.
Moreover, if $A$ is an algebraic structure, $\sigma(i)\in {\bf Con}(A)$, $i\in I$ and $\mathcal{D}=\{ K_j| j\in J\}$ then the restriction of $\Phi (\sigma /\mathcal{D}^*)$ to $A$ equals $\vee _{j\in J}(\wedge _{k\in K_j}\sigma (k))$.
\end{abstract}

\section{Introduction}

Let $I$ be a nonempty set and ${{\D}}$ an ultrafilter over $I$, that is, ${\D}$ is a subset of the power set $P(I)$ of $I$ with the following properties:
\begin{itemize}
\item[(1)] $I\in {\D}$ but $\emptyset \notin {\D}$;
\item[(2)] $A\cap B\in {\D}$ for any $A, B \in {\D}$;
\item[(3)] If $A\in {\D}$ and $C\subseteq I$ with $A\subseteq C$ then $C\in {\D}$;
\item[(4)] For each $A\subseteq I$, $A\in \mathcal{D}$ or $I\setminus A\in \mathcal{D}$.
\end{itemize}

\medskip
We note that $(1)$ and $(2)$ together imply that, for each $A\subseteq I$, exactly one of the sets $A$ , $I\setminus A$ belongs to $\mathcal{D}$.
Moreover, by Corollary 3.13 of Chapter IV of \cite{Burris}, condition $(4)$ can be replaced by the following condition:

\begin{itemize}
\item[($4^*$)] For any $A, B\subseteq I$, the assumption $A\cup B\in {\D}$ implies $A\in {\D}$ or $B\in {\D}$.
\end{itemize}

\medskip

Let $A_i=(A_i; \Omega)$, $i\in I$ be arbitrary similar algebraic structures and $\sigma _i$ be a congruence on $A_i$, $i\in I$. Let
$\Pi (A_i|i\in I)$ denote the direct product of $A_i$.
Let ${\bf Con}(A_i)$ denote the $\wedge$-semilattice of all congruences of $A_i$, $i\in I$. Let $\sigma \in \Pi ({\bf Con}(A_i)|i\in I)$ be an element for which $\sigma (i)=\sigma _i$, $i\in I$. Let $\Pi (\sigma (i)|i\in I)$ denote the relation on the direct product $\Pi (A_i|i\in I)$ defined by
$(a, b)\in \Pi (\sigma (i)|i\in I)$ if and only if $\{ i\in I| (a(i), b(i))\in \sigma (i)\}\in \mathcal{D}$. By the dual of Lemma 3 of \S 8 of \cite{Gratzer}, it is easy to see that $\Pi (\sigma (i)|i\in I)$ is a congruence on $\Pi (A_i|i\in I)$. Recall that this congruence is the ultraproduct congruence (denoted by ${\D}^*$) on
$\Pi (A_i|i\in I)$ if $\sigma (i)$ is the identity relation on $A_i$ for all $i\in I$. It is evident, that ${\D}^*\subseteq \Pi (\sigma (i)|i\in I)$ for arbitrary congruences $\sigma (i)$ on $A_i$, $i\in I$. Then we can consider the congruence
$\Pi (\sigma (i)|i\in I)/{\D}^*$ on the ultraproduct $\Pi _{\D}(A_i|i\in I)=\Pi (A_i|i\in I)/{\D}^*$ (see Definition 6.13 and Lemma 6.14 of Chapter II of \cite{Burris}). This congruence will be denoted by $\Pi _{\mathcal{D}}(\sigma (i)|i\in I)$.
In this paper we examine congruences $\Pi _{\mathcal{D}}(\sigma (i)|i\in I)$ on the ultraproduct $\Pi _{\D}(A_i|i\in I)$.

Let the ultraproduct congruence on $\Pi ({\bf Con}(A_i)|i\in I)$ denoted by also $\mathcal{D}^*$.
In Section 2, we show that $\Phi : \sigma /\mathcal{D}^* \mapsto \Pi _{\mathcal{D}}(\sigma(i)|i\in I)$ is an embedding of the $\wedge$-semilattice $\Pi _{\mathcal{D}}({\bf Con}(A_i)|i\in I)$ into the $\wedge$-semilattice ${\bf Con}(\Pi _{\D}(A_i|i\in I)$, where $\sigma /\mathcal{D}^*$ denotes the ultraproduct congruence class on the direct product $\Pi ({\bf Con}(A_i)|i\in I)$ containing $\sigma=(\sigma(i))_{i\in I}$. In Section 3, we prove that the factor algebra $\Pi _{\mathcal{D}}(A_i|i\in I)/\Pi _{\mathcal{D}}(\sigma (i)|i\in I)$ is isomorphic to the ultraproduct $\Pi _{\mathcal{D}}(A_i/\sigma(i)|i\in I)$. In Section 4, we apply our results for the ultrapower $\Pi _{\mathcal{D}}(A|i\in I)$ of an algebraic structure $A$. Since $A$ can be embedded into its ultrapower $\Pi _{\mathcal{D}}(A|i\in I)$ then, for an arbitrary family $\{\sigma (i)|i\in I\}$ of congruences on $A$, we can consider the restriction $\Pi _{\mathcal{D}}(\sigma (i)| i\in I)|A$ of the congruence $\Pi _{\mathcal{D}}(\sigma (i)| i\in I)$ on
$\Pi _{\mathcal{D}}(A|i\in I)$ to $A$. We show that if $\{ \sigma (i)| i\in I\}$ is an arbitrary family of congruences on $A$ and
${\mathcal{D}}=\{ K_j|j\in J\}$ then $\cup _{j\in J}(\cap _{k\in K_j}\sigma _k)$ is a congruence on $A$ such that $\Pi _{\mathcal{D}}(\sigma (i)| i\in I)|A=\cup _{j\in J}(\cap _{k\in K_j}\sigma _k)=\vee _{j\in J}(\wedge _{k\in K_j}\sigma (k))$.

\section{An embedding theorem}


\begin{theorem}\label{th1} Let $I$ be a non-empty set and $\mathcal{D}$ an ultrafilter over $I$. Then the mapping $\Phi : \sigma/\mathcal{D}^* \mapsto \Pi _{\D}(\sigma (i)|i\in I)$ is an embedding of the $\wedge$-semilattice $\Pi _{\D}({\bf Con}(A_i)|i\in I)$ into the $\wedge$-semilattice ${\bf Con}(\Pi _{\D}(A_i|i\in I))$ for arbitrary similar algebraic structures $A_i$, $i\in I$.
\end{theorem}

\noindent
{\bf Proof}.
Let $\alpha =(\alpha (i))_{i\in I}$ and $\beta =(\beta (i))_{i\in I}$ be arbitrary elements of the direct product $\Pi ({\bf Con}(A_i)|i\in I)$.
First we show that
\[\alpha /\mathcal{D}^*=\beta /\mathcal{D}^* \Longleftrightarrow
\Pi  _{\D}(\alpha(i)|i\in I)=\Pi  _{\D}(\beta(i)|i\in I).\]
By the Correspondence Theorem (Theorem 6.20 of Chapter II of \cite{Burris}), it is sufficient to show that
\[\alpha /\mathcal{D}^*=\beta /\mathcal{D}^* \Longleftrightarrow
\Pi (\alpha(i)|i\in I)=\Pi (\beta(i)|i\in I).\]
Assume $\alpha /\mathcal{D}^*=\beta /\mathcal{D}^*$. Then \[A=\{i\in I: \alpha (i)=\beta (i)\}\in {\D}.\]
If \[(a, b)\in \Pi (\alpha(i)|i\in I)\] for some elements $a=(a(i))_{i\in I}$ and $b=(b(i))_{i\in I}$ of $\Pi (A_i|i\in I)$ then
\[B=\{i\in I: (a(i), b(i))\in \alpha (i)\}\in {\D}\] and so, for every $i\in A\cap B\in {\D}$,
\[(a(i), b(i))\in \beta (i)\] from which it follows that \[\{i\in I: (a(i), b(i))\in \beta (i)\}\in {\D}.\] Thus
\[(a, b)\in \Pi (\beta(i)|i\in I).\] Hence
\[\Pi (\alpha(i)|i\in I)\subseteq \Pi (\beta(i)|i\in I).\] We can prove
\[\Pi (\beta(i)|i\in I)\subseteq \Pi (\alpha(i)|i\in I)\] in a similar way. Thus
\[\Pi (\alpha(i)|i\in I)=\Pi (\beta(i)|i\in I).\]

 To prove the converse, assume
 \[\Pi (\alpha(i)|i\in I)=\Pi (\beta(i)|i\in I).\] We show that \[K=\{i\in I| \alpha(i)=\beta (i)\}\in {\D}.\]
Assume, in an indirect way, that $K\notin \mathcal{D}$. Then $I\setminus K\in \mathcal{D}$.

First assume \[R=\{i\in I: \alpha (i)\subset \beta (i)\}\in {\D}.\] Then $(I\setminus R)\notin {\D}$.  For $i\in R$, let $(a_i, b_i)$ be an element of $A_i\times A_i$ such that $(a_i, b_i)\in \beta (i)$ and $(a_i, b_i)\notin \alpha (i)$. For the indexes $i\in I\setminus R$, let $(a_i, b_i)\in A_i\times A_i$ be arbitrary. Let $a$ and $b$ be the elements of $\Pi (A_i|i\in I)$ for which $a(i)=a_i$ and $b(i)=b_i$. As
$(a(i), b(i))\in \beta (i)$ for every $i\in R\in {\D}$, we have \[\{i\in I: (a(i), b(i))\in \beta (i)\}\in {\D}\]
and so \[(a, b)\in \Pi  (\beta(i)|i\in I).\] We show that $(a, b)\notin \Pi (\alpha(i)|i\in I)$.
Assume, in an indirect way, that \[(a, b)\in \Pi  (\alpha(i)|i\in I).\] Then
\[V=\{i\in I: (a(i), b(i))\in \alpha (i)\}\in {\D}.\] As $V\subseteq I\setminus R$, we have $I\setminus R\in {\D}$ and so
$R\notin {\D}$ which is a contradiction. Hence \[(a, b)\notin \Pi (\alpha(i)|i\in I).\] This implies
\[\Pi (\alpha(i)|i\in I)\neq \Pi \beta(i)|i\in I)\] which is a contradiction.

Next consider the case when
\[R=\{i\in I: \alpha (i)\subset \beta (i)\}\notin {\D}.\]
Then
\[\{i\in I: \alpha (i)\subseteq \beta (i)\}=K\cup R\notin {\D}.\]
For every $i\in I\setminus (K\cup R)\in \mathcal{D}$, there are elements $a_i, b_i$ of $A_i$ such that
\[(a_i, b_i)\in \alpha (i),\quad (a_i, b_i)\notin \beta (i).\] For an index $i\in K\cup R$, let the elements $a_i, b_i\in A_i$ be arbitrary.
Let $a$ and $b$ be the elements of $\Pi (A_i|i\in I)$ for which $a(i)=a_i$ and $b(i)=b_i$.
It is clear that
\[(a, b)\in \Pi (\alpha(i)|i\in I).\] We show that $(a, b)\notin \Pi  _{\D}(\beta(i)|i\in I)$.
Assume, in an indirect way, that
\[(a, b)\in \Pi (\beta(i)|i\in I).\] Then \[\{i\in I: (a(i), b(i))\in \beta (i)\}\in {\D}.\] As
\[\{i\in I: (a(i), b(i))\in \beta (i)\}\subseteq K\cup R,\] we have $K\cup R\in {\D}$ which is a contradiction. Thus
\[(a, b)\notin \Pi (\beta(i)|i\in I).\] Hence
\[\Pi (\alpha(i)|i\in I)\neq \Pi (\beta(i)|i\in I)\] which is a contradiction. We got a contradiction in both cases. Thus
\[K=\{i\in I: \alpha (i)=\beta (i)\}\in {\D}.\]  Consequently
$\alpha /\mathcal{D}^*=\beta /\mathcal{D}^*$.

Thus $\Phi : \sigma/\mathcal{D}^* \mapsto \Pi _{\D}(\sigma (i)|i\in I)$ is a well defined injective mapping.
We show that $\Phi$ is a homomorphism. For arbitrary $a, b\in \Pi (A_i|i\in I)$, let
\[I_{\alpha \wedge \beta}=\{ i\in I: (a(i), b(i))\in (\alpha\wedge \beta)(i)=\alpha (i)\wedge \beta (i)\},\]
\[ I_{\alpha}=\{ i\in I: (a(i), b(i))\in \alpha (i)\}\]
and
\[ I_{\beta}=\{ i\in I: (a(i), b(i))\in \beta (i)\}.\]
It is easy to see that \[I_{\alpha \wedge \beta}=I_{\alpha}\cap I_{\beta}.\]
As $\mathcal{D}$ is an ultrafilter, $I_{\alpha}\cap I_{\beta}\in {\D}$ if and only if $I_{\alpha}\in {\D}$ and $I_{\beta }\in {\D}$. From the above it follows that

\[(a, b)\in \Pi  ((\alpha \wedge \beta)(i)|i\in I)\Longleftrightarrow (a, b)\in \Pi  (\alpha (i)|i\in I)\wedge \Pi  (\beta(i)|i\in I).\]
Hence \[\Pi  ((\alpha \wedge \beta)(i)|i\in I)=\Pi (\alpha (i)|i\in I)\wedge \Pi (\beta(i)|i\in I).\]
By the Correspondence Theorem (Theorem 6.20 of Chapter II of \cite{Burris}), we have
\[\Pi _{\mathcal{D}}((\alpha \wedge \beta)(i)|i\in I)=\Pi _{\mathcal{D}}(\alpha (i)|i\in I)\wedge \Pi  _{\D}(\beta(i)|i\in I).\]
Then
\[\Phi(\alpha /\mathcal{D}^*\wedge \beta /\mathcal{D}^*)=\]
\[\Phi ((\alpha \wedge \beta)/\mathcal{D}^*)=\Pi _{\mathcal{D}}((\alpha \wedge \beta)(i)|i\in I)=
\Pi _{\mathcal{D}}(\alpha(i)|i\in I)\wedge \Pi _{\mathcal{D}} (\beta(i)|i\in I)=\]
\[=\Phi (\alpha/\mathcal{D}^*)\wedge \Phi (\beta/\mathcal{D}^*),\]
and so $\Phi$ is a homomorphism. Thus $\Phi$ is an embedding of the $\wedge$-semilattice $\Pi  _{{\D}}({\bf Con}(A_i)|i\in I)$ into the $\wedge$-semilattice
${\bf Con}(\Pi  _{{\D}}(A_i|i\in I))$.\hfill\openbox

\section{An isomorphism theorem}

\begin{theorem}\label{th2} Let $I$ be a non-empty set and $\mathcal{D}$ an ultrafilter over $I$. Then, for arbitrary similar algebraic structures $A_i$, $i\in I$ and an arbitrary congruence $\sigma (i)$ on $A_i$ ($i\in I$), the factor algebra $\Pi _{\mathcal{D}}(A_i|i\in I)/\Pi _{\mathcal{D}}(\sigma (i)|i\in I)$ is isomorphic to the ultraproduct $\Pi _{\D}(A_i/\sigma (i)|i\in I)$ of the factor algebras $A_i/\sigma (i)$.
\end{theorem}

\noindent
{\bf Proof}.
As $\Pi _{\mathcal{D}}(A_i|i\in I)/\Pi _{\mathcal{D}}(\sigma (i)|i\in I)\cong \Pi (A_i|i\in I)/\Pi (\sigma (i)|i\in I)$ by the Second Isomorphic Theorem (Theorem 6.15 of Chapter II of \cite{Burris}), it is sufficient to show that the factor algebra $\Pi (A_i|i\in I)/\Pi (\sigma (i)|i\in I)$ is isomorphic to the ultraproduct
$\Pi _{\D}(A_i/\sigma (i)|i\in I)$.
For an element $a=(a(i))_{i\in I}$ of the direct product $\Pi (A_i|i\in I)$, let
\[\Delta (a)=([a(i)]_{\sigma (i)})_{i\in I}/\mathcal{D}^*,\] where $[a(i)]_{\sigma (i)}$ denotes the $\sigma (i)$-class of $A_i$ containing $a(i)$. It is obvious that $\Delta: \Pi (A_i|i\in I) \mapsto \Pi _{\D}(A_i/\sigma (i)|i\in I)$ is a well defined surjective mapping. Let $a_j=(a_j(i))_{i\in I}$ ($j=1, \dots , n$) be arbitrary elements of $\Pi (A_i|i\in I)$. If $\omega _n\in \Omega$ is an $n$-ary operation then, for arbitrary $i\in I$,
\[\omega _n(\Delta (a_1), \dots , \Delta (a_n))=\omega _n(([a_1(i)]_{\sigma (i)})_{i\in I}/{{\D}^*}, \dots ,
([a_n(i)]_{\sigma (i)})_{i\in I}/{{\D}^*})=\]
\[=\omega _n(([a_1(i)]_{\sigma (i)})_{i\in I}, \dots , ([a_n(i)]_{\sigma (i)})_{i\in I})/{{\D}^*}=\]
\[=(\omega _n([a_1(i)]_{\sigma (i)}, \dots , [a_n(i)]_{\sigma (i)}))_{i\in I}/{{\D}^*}=
([\omega _n(a_1(i), \dots , a_n(i))] _{\sigma (i)})_{i\in I}/{{\D}^*}=\]
\[=(([\omega _n(a_1, \dots , a_n)(i)] _{\sigma (i)})_{i\in I})/{{\D}^*}=\Delta(\omega _n(a_1, \dots , a_n)).\]
Thus $\Delta$ is a homomorphism. For elements $a=(a(i))_{i\in I}$ and $b=(b(i))_{i\in I}$ of $\Pi (A_i|i\in I)$,
\[\Delta (a)=\Delta (b)\] if and only if \[([a(i)]_{\sigma (i)})_{i\in I}/{{\D}^*}=([b(i)]_{\sigma (i)})_{i\in I}/{{\D}^*},\] that is,
\[\{ i\in I| [a(i)]_{\sigma (i)}=[b(i)]_{\sigma (i)}\}\in \mathcal{D}.\] This last condition is equivalent to the condition that
\[\{ i\in I| (a(i), b(i))\in \sigma (i)\}\in \mathcal{D},\] that is, \[(a, b)\in \Pi (\sigma (i)|i\in I).\] Thus the kernel of $\Delta$ is
$\Pi (\sigma (i)|i\in I)$. Our assertion follows from the Homomorphism Theorem (Theorem 6.12 of Chapter II of \cite{Burris}).
\hfill\openbox

\section{Ultrapowers}

Let $A$ be an algebraic structure. If $\mathcal{D}$ is an ultrafilter over a non-empty set $I$, we can consider the ultrapower of $A$ modulo $\mathcal{D}$ as the ultraproduct $\Pi _{\mathcal {D}}(A_i|i\in I)$, where $A_i=A$ for all $i\in I$.  For an arbitrary element $a\in A$, let
$\xi (a)$ denote the $\mathcal{D}^*$-class of the direct product $\Pi (A|i\in I)$ which contains the constant function with value $a$.
It is known (Lemma 2.10 of Chapter V of \cite{Burris}) that $\xi :a\mapsto \xi (a)$ is an embedding (the natural embedding) of $A$ into the ultrapower $\Pi _{\mathcal{D}}(A|i\in I)$. Identify $A$ and $\xi (A)$.

Let $\sigma (i)$, $i\in I$ be arbitrary congruences on $A$. Let $\Pi _{\mathcal{D}}(\sigma (i)|i\in I)|A$
denote the restriction of the congruence $\Pi _{\mathcal{D}}(\sigma (i)|i\in I)$ to $A$.
It is clear that, for some $a, b\in A$, $(a, b)\in \Pi _{\mathcal{D}}(\sigma (i)|i\in I)|A$ if and only if $\{ i\in I| (a, b)\in \sigma (i)\}\in {\mathcal{D}}$.

\medskip

\begin{theorem}\label{equiv} Let $I$ be a non-empty set and $\mathcal{D}=\{K_j|j\in J\}$ an ultrafilter over $I$.
Then, for an arbitrary family $\{ \sigma (i)| i\in I\}$ of congruences $\sigma (i)$ on an algebraic structure $A$, $\cup _{j\in J}(\cap _{k\in K_j}\sigma (k))$ is a congruence on $A$ such that \\
$\Pi _{\mathcal{D}}(\sigma (i)|i\in I)|_A=\cup _{j\in J}(\cap _{k\in K_j}\sigma (k))=\vee _{j\in J}(\wedge _{k\in K_j}\sigma (k))$.
\end{theorem}

\noindent
{\bf Proof}. For a couple $(a, b)$ of elements $a$ and $b$ of $A\subseteq \Pi _{\mathcal {D}}(A|i\in I)$, let
\[K_{a, b}=\{i\in I| (a, b)\in \sigma (i)\}.\]
Assume $(a, b)\in \Pi _{\mathcal{D}}(\sigma (i)|i\in I)$ for some $a, b\in A$. Then $K_{a, b}\in {\mathcal{D}}$
and so \[(a, b)\in \cap _{k\in K_{a, b}}\sigma(k)\] from which it follows that
\[(a, b)\in \cup _{j\in J}(\cap _{k\in K_j}\sigma (k)).\] Thus
\[\Pi _{\mathcal{D}}(\sigma (i)|i\in I)|_A\subseteq \cup _{j\in J}(\cap _{k\in K_j}\sigma (k)).\]

Conversely, assume \[(a, b)\in \cup _{j\in J}(\cap _{k\in K_j}\sigma (k))\] for some $a, b\in A$. Then there is an index $j\in J$ such that
\[(a, b)\in (\cap _{k\in K_j}\sigma (k)).\]
Thus $K_j\subseteq K_{a, b}$ and so $K_{a, b}\in {\mathcal{D}}$.
Hence \[(a, b)\in \Pi _{\mathcal{D}}(\sigma (i)|i\in I)|_A.\] Thus
\[\cup _{j\in J}(\cap _{k\in K_j}\sigma (k)\subseteq \Pi _{\mathcal{D}}(\sigma (i)|i\in I)|_A.\]
Consequently
\[\Pi _{\mathcal{D}}(\sigma (i)|i\in I)|_A=\cup _{j\in J}(\cap _{k\in K_j}\sigma (k)).\]
Thus $\cup _{j\in J}(\cap _{k\in K_j}\sigma (k))$ is a congruence on $A$. It is obvious that
$\cup _{j\in J}(\cap _{k\in K_j}\sigma (k))$ is a common upper bound of the congruences $\sigma _j=\cap _{k\in K_j}\sigma (k)$, $j\in J$. Moreover, for every common upper bound $\beta$ of $\sigma _j$, $j\in J$, we have $\cup _{j\in J}(\cap _{k\in K_j}\sigma (k))\subseteq \beta$ and so
$\cup _{j\in J}(\cap _{k\in K_j}\sigma (k))=\vee _{j\in J}(\wedge _{k\in K_j}\sigma (k))$.
Hence
\[\Pi _{\mathcal{D}}(\sigma (i)|i\in I)|_A=\cup _{j\in J}(\cap _{k\in K_j}\sigma (k))=\vee _{j\in J}(\wedge _{k\in K_j}\sigma (k)).\]
\hfill\openbox

\medskip
\noindent
Department of Algebra

\noindent
Budapest University of Technology and Economics

\noindent
P.O. Box 91, 1521 Budapest, Hungary

\noindent
e-mail: nagyat@math.bme.hu

\end{document}